\documentclass[12pt]{amsart}
\usepackage{amsmath,amssymb,amsfonts,amscd}
\theoremstyle{plain}
\newtheorem{thm}{Theorem}
\newtheorem{lem}[thm]{Lemma}

\newtheorem{prop}[thm]{Proposition}

\numberwithin{thm}{section}
\numberwithin{equation}{section}

\newcommand{\eq}[2]{\begin{equation}\label{#1}#2 \end{equation}}
\newcommand{\ml}[2]{\begin{multline}\label{#1}#2 \end{multline}}
\newcommand{\ga}[2]{\begin{gather}\label{#1}#2 \end{gather}}

\newcommand{\inj}{\hookrightarrow}

\newcommand{\Spec}{{\rm Spec \,}}

\newcommand{\sO}{{\mathcal O}}

\newcommand{\A}{{\Bbb A}}

\newcommand{\G}{{\Bbb G}}

\renewcommand{\P}{{\Bbb P}}
\newcommand{\Q}{{\Bbb Q}}

\newcommand{\T}{{\Bbb T}}

\newcommand{\X}{{\Bbb X}}

\newcommand{\Z}{{\Bbb Z}}

\newcommand{\ve}{{\varepsilon}}

\title{Motives Associated to Sums of Graphs}
\author{Spencer Bloch}\address{Dept. of Mathematics, University of Chicago, Chicago, IL 60637,
USA\\
E-mail address: bloch@math.uchicago.edu}
\begin{document}
\maketitle
\section{Introduction}In quantum field theory, the path integral is interpreted perturbatively as a sum indexed by graphs. The coefficient (Feynman amplitude) associated to a graph $\Gamma$ is a period associated to the motive given by the complement of a certain hypersurface $X_\Gamma$ in projective space. Based on considerable numerical evidence, Broadhurst and Kreimer suggested \cite{BK} that the Feynman amplitudes should be sums of multi-zeta numbers. On the other hand, Belkale and Brosnan \cite{BB} showed that the motives of the $X_\Gamma$ were not in general mixed Tate. 

A recent paper of Aluffi and Marcolli \cite{AM} studied the images $[X_\Gamma]$ of graph hypersurfaces in the Grothendieck ring $K_0(Var_k)$ of varieties over a field $k$. Let $\Z[\A^1_k] \subset K_0(Var_k)$ be the subring generated by $1= [\Spec k]$ and $[\A^1_k]$. It follows from \cite{BB} that $[X_\Gamma]\not\in \Z[\A^1_k]$ for many graphs $\Gamma$. 

Let $n\ge 3$ be an integer. In this note we consider a sum $S_n\in K_0(Var_k)$ of $[X_\Gamma]$ over all connected graphs $\Gamma$ with $n$ vertices, no multiple edges, and no tadpoles (edges with just one vertex). (There are some subtleties here. Each graph $\Gamma$ appears with multiplicity $n!/|Aut(\Gamma)|$. For a precise definition of $S_n$ see \eqref{5.1} below.) Our main result is
\begin{thm}$S_n \in \Z[\A^1_k]$. 
\end{thm}

For applications to physics, one would like a formula for sums over all graphs with a given loop order. I do not know if such a formula could be proven by these methods. 

Dirk Kreimer explained to me the physical interest in considering sums of graph motives, and I learned about $K_0(Var_k)$ from correspondence with H. Esnault. Finally, the recently paper of Aluffi and Marcolli \cite{AM} provides a nice exposition of the general program. 

\section{Basic Definitions} Let $E$ be a finite set, and let 
\eq{2.1}{0 \to H \to \Q^E \to W \to 0;\quad 0 \to W^\vee \to \Q^E \to H^\vee \to 0
}
be dual exact sequences of vector spaces. For $e\in E$, let $e^\vee: \Q^E \to \Q$ be the dual functional, and let $(e^\vee)^2$ be the square, viewed as a quadratic function. By restriction, we can view this as a quadratic function either on $H$ or on $W^\vee$. Choosing bases, we get symmetric matrices $M_e$ and $N_e$. Let $A_e, e\in E$ be variables, and consider the homogeneous polynomials
\eq{2.2}{\Psi(A) = \det(\sum A_eM_e);\quad \Psi^\vee(A) = \det(\sum A_eN_e).
} 
\begin{lem}\label{lem2.1} $\Psi(\ldots A_e,\ldots) = c\prod_{e\in E} A_e \Psi^\vee(\ldots A_e^{-1},\ldots)$, where $c\in k^\times$. 
\end{lem}
\begin{proof}This is proposition 1.6 in \cite{BEK}. 
\end{proof}

Let $\Gamma$ be a graph. Write $E, V$ for the edges and vertices of $\Gamma$. We have an exact sequence
\eq{}{0 \to H_1(\Gamma,\Q) \to \Q^E \xrightarrow{\partial} \Q^V \to H_0(\Gamma,\Q) \to 0.
}
We take $H = H_1(\Gamma)$ and $W=\text{Image}(\partial)$ in \eqref{2.1}. The resulting polynomials $\Psi = \Psi_\Gamma,\ \Psi^\vee = \Psi_\Gamma^\vee$ as in \eqref{2.2} are given by \cite{BEK}
\eq{}{\Psi_\Gamma = \sum_{t\in T} \prod_{e\not\in t} A_e;\quad \Psi_\Gamma^\vee = \sum_{t\in T} \prod_{e\in t} A_e.
}
Here $T$ is the set of {\it spanning trees} in $\Gamma$. 

\begin{lem}\label{lem2.2} Let $e \in \Gamma$ be an edge. Let $\Gamma/e$ be the graph obtained from $\Gamma$ by shrinking $e$ to a point and identifying the two vertices. We do not consider $\Gamma/e$ in the degenerate case when $e$ is a loop, i.e. if the two vertices coincide. Let $\Gamma-e$ be the graph obtained from $\Gamma$ by cutting $e$. We do not consider $\Gamma-e$ in the degenerate case when cutting $e$ disconnects $\Gamma$ or leaves an isolated vertex. Then 
\ga{2.5}{\Psi_{\Gamma/e} = \Psi_\Gamma|_{A_e=0};\quad \Psi_{\Gamma-e} = \frac{\partial}{\partial A_e}\Psi_\Gamma. \\
\Psi_{\Gamma/e}^\vee =\frac{\partial}{\partial A_e}\Psi_\Gamma^\vee; \quad \Psi_{\Gamma-e}^\vee = \Psi_\Gamma^\vee |_{A_e=0}. \label{2.6}
}
(In the degenerate cases, the polynomials on the right in \eqref{2.5} and \eqref{2.6} are zero.) 
\end{lem}
\begin{proof} The formulas in \eqref{2.5} are standard \cite{BEK}. The formulas \eqref{2.6} follow easily using lemma \ref{lem2.1}. (In the case of graphs, the constant $c$ in the lemma is $1$.) 
\end{proof}

More generally, we can consider strings of edges $e_1,\dotsc,e_p \in \Gamma$. If at every stage we have a nondegenerate situation we can conclude inductively
\eq{2.7}{\Psi_{\Gamma-e_1-\cdots -e_p}^\vee = \Psi_\Gamma^\vee|_{A_{e_1}=\cdots = A_{e_p}=0}
}
In the degenerate situation, the polynomial on the right will vanish, i.e. $X_\Gamma$ will contain the linear space $A_{e_1}=\cdots =A_{e_p}=0$. 

For example, let $\Gamma = e_1\cup e_2\cup e_3$ be a triangle, with one loop and three vertices. We get the following polynomials
\ga{2.8}{\Psi_\Gamma = A_{e_1}+A_{e_2}+A_{e_3}; \quad \Psi_\Gamma^\vee = A_{e_1}A_{e_2}+A_{e_2}A_{e_3}+A_{e_1}A_{e_3} \\
\Psi_{\Gamma-e_i} = 1;\quad \Psi_{\Gamma-e_i}^\vee = A_{e_j}A_{e_k} = \Psi_\Gamma^\vee|_{A_{e_i}=0} \label{2.9}
}
The sets $\{e_i,e_j\}$ are degenerate because cutting two edges will leave an isolated vertex. 

\section{The Grothendieck Group and Duality}

Recall $K_0(Var_k)$ is the free abelian group on generators isomorphism classes $[X]$ of quasi-projective $k$-varieties and relations 
\eq{1.1}{[X] = [U]+[Y]; \quad U \stackrel{\text{open}}{\inj} X,\ Y= X-U.
}
In fact, $K_0(Var_k)$ is a commutative ring with multiplication given by cartesian product of $k$-varieties. Let $\Z[\A^1_k] \subset K_0(Var_k)$ be the subring generated by $1= [\Spec k]$ and $[\A^1_k]$. 
Let $\P_\Gamma$ be the projective space with homogeneous coordinates $A_e, e\in E$. We write $X_\Gamma: \Psi_\Gamma=0,\ X_\Gamma^\vee : \Psi_\Gamma^\vee = 0$ for the corresponding hypersurfaces in $\P_\Gamma$. We are interested in the classes $[X_\Gamma], [X_\Gamma^\vee] \in K_0(Var_k)$. 

Let $\Delta: \prod_{e\in E} A_e=0$ in $\P_\Gamma$, and let $\T=\T_\Gamma = \P_\Gamma-\Delta$ be the torus. Define
\eq{}{X_\Gamma^0 = X_\Gamma\cap \T_\Gamma;\quad X_\Gamma^{\vee,0} = X_\Gamma^\vee\cap \T_\Gamma.
} 
Lemma \ref{lem2.1} translates into an isomorphism (Cremona transformation) 
\eq{3.3}{X_\Gamma^0 \cong  X_\Gamma^{\vee,0}.
}
(In fact, this is valid more generally for the setup of \eqref{2.1} and \eqref{2.2}.) We can stratify $X_\Gamma^{\vee}$ by intersecting with the toric stratification of $\P_\Gamma$ and write
\eq{3.4}{[X_\Gamma^{\vee}] = \sum_{\{e_1,\dotsc,e_p\} \subset E}[(X_\Gamma^{\vee}\cap \{A_{e_1}=\cdots = A_{e_p}=0\})^0] \in K_0(Var_k)
}
where the sum is over all subsets of $E$, and superscript $0$ means the open torus orbit where $A_e \neq 0, e\not\in \{e_1,\dotsc,e_p\}$. We call a subset $\{e_1,\dotsc,e_p\} \subset E$ degenerate if $\{A_{e_1}=\cdots = A_{e_p}=0\} \subset X_\Gamma^\vee$. Since $[\G_m] = [\A^1]-[pt] \in K_0(Var_k)$ we can rewrite \eqref{3.4}
\eq{3.5}{[X_\Gamma^{\vee}] = \sum_{\substack{\{e_1,\dotsc,e_p\} \subset E \\ \text{nondegenerate}}} [(X_\Gamma^{\vee}\cap \{A_{e_1}=\cdots = A_{e_p}=0\})^0] + t
}
where $t \in \Z[\A^1] \subset K_0(Var_k)$. Now using \eqref{2.7} and \eqref{3.3} we conclude 
\eq{3.6}{[X_\Gamma^{\vee}] = \sum_{\substack{\{e_1,\dotsc,e_p\} \subset E \\ \text{nondegenerate}}} [(X_{\Gamma-\{e_1,\dotsc,e_p\}}^0] + t.
}

\section{Complete Graphs}

Let $\Gamma_n$ be the complete graph with $n\ge 3$ vertices. Vertices of $\Gamma_n$ are written $(j),\ 1\le j\le n$, and edges $e_{ij}$ with $1\le i<j\le n$. We have $\partial e_{ij} = (j)-(i)$. 
\begin{prop}\label{prop4.1} We have $[X_{\Gamma_n}^\vee] \in \Z[\A^1_k]$. 
\end{prop}
\begin{proof}Let $\Q^{n,0} \subset \Q^n$ be row vectors with entries which sum to $0$. We have
\eq{}{0 \to H_1(\Gamma_n) \to \Q^E \xrightarrow{\partial} \Q^{n,0} \to 0.
}
In a natural way, $(\Q^{n,0})^\vee = \Q^n/\Q$. Take as basis of $\Q^n/\Q$ the elements $(1),\dotsc,(n-1)$. As usual, we interpret the $(e_{ij}^\vee)^2$ as quadratic functions on $\Q^n/\Q$. We write $N_e$ for the corresponding symmetric matrix. 
\begin{lem}The $N_{e_{ij}}$ form a basis for the space of all $(n-1)\times (n-1)$ symmetric matrices. 
\end{lem}
\begin{proof}[Proof of lemma] The dual map $\Q^n/\Q \to \Q^E$ carries 
\eq{}{(k) \mapsto \sum_{\mu>k} -e_{k\mu} + \sum_{\nu<k} e_{\nu k};\quad k\le n-1. 
}We have
\eq{}{(e_{ij}^\vee)^2(\sum_{k=1}^{n-1}a_k\cdot (k)) = \begin{cases}a_i^2-2a_ia_j+a_j^2 & i<j<n \\
a_i^2 & j=n. \end{cases}
}
It follows that if $j<n$, $N_{e_{ij}}$ has $-1$ in positions $(ij)$ and $(ji)$ and $+1$ in positions $(ii), (jj)$ (resp. $N_{in}$ has $1$ in position $(ii)$ and zeroes elsewhere). These form a basis for the symmetric $(n-1)\times (n-1)$ matrices. 
\end{proof}

It follows from the lemma that $X_{\Gamma_n}^\vee$ is identified with the projectivized space of $(n-1)\times (n-1)$ matrices of rank $\le n-2$. In order to compute the class in the Grothendieck group we detour momentarily into classical algebraic geometry. For a finite dimensional $k$-vector space $U$, let $\P(U)$ be the variety whose $k$-points are the lines in $U$. For a $k$-algebra $R$, the $R$-points $\Spec R \to \P(U)$ are given by pairs $(L, \phi)$ where $L$ on $Spec\ R$ is a line bundle and $\phi: L\inj U \otimes_k R$ is a locally split embedding. 

Suppoose now $U = Hom(V,W)$. We can stratify $\P(Hom(V,W))= \coprod_{p> 0}\P(Hom(V,W))^p$ according to the rank of the homomorphism. Looking at determinants of minors makes it clear that $\P(Hom(V,W))^{\le p}$ is closed. Let $R$ be a local ring which is a localization of a $k$-algebra of finite type, and let $a$ be an $R$-point of $\P(Hom(V,W))^p$. Choosing a lifting $b$ of the projective point $a$, we have
\eq{2.12}{0 \to \ker(b) \to V\otimes R \xrightarrow{b} W\otimes R \to \text{coker}(b) \to 0,
}
and $\text{coker}(b)$ is a finitely generated $R$-module of constant rank $\dim W - p$ which is therefore necessarily free. 

Let $Gr(\dim V-p,V)$ and $Gr(p,W)$ denote the Grassmann varieties of subspaces of the indicated dimension in $V$ (resp. $W$). On $Gr(\dim V-p,V)\times Gr(p,W)$ we have rank $p$ bundles $E,F$ given respectively by the pullbacks of the universal quotient on $Gr(\dim V-p,V)$ and the universal subbundle on $Gr(p,W)$. It follows from the above discussion that
\eq{}{\P(Hom(V,W))^p = \P(\text{Isom}(E,F)) \subset \P(Hom(E,F)).
}

Suppose now that $W = V^\vee$. Write $\langle\ ,\ \rangle: V \otimes V^\vee \to k$ for the canonical bilinear form. We can identify $Hom(V,V^\vee)$ with bilinear forms on $V$ 
\eq{}{\rho: V \to V^\vee \leftrightarrow (v_1,v_2) \mapsto \langle v_1,\rho(v_2)\rangle. 
}
Let $SHom(V,V^\vee) \subset Hom(V,V^\vee)$ be the subspace of $\rho$ such that the corresponding bilinear form on $V$ is symmetric. Equivalently, $Hom(V,V^\vee) = V^{\vee,\otimes 2}$ and $SHom(V,V^\vee) = Sym^2(V^\vee) \subset V^{\vee,\otimes 2}$. 

For $\rho$ symmetric as above, one seees easily that $\rho(V) = \ker(V)^\perp$ so there is a factorization
\eq{2.15}{V \to V/\ker(\rho) \xrightarrow{\cong}(V/\ker(\rho))^\vee  =  \ker(\rho)^\perp \inj V^\vee. 
}
The isomorphism in \eqref{2.15} is also symmetric. 

Fix an identification $V=k^{n}$ and hence $V=V^\vee$. A symmetric map is then given by a symmetric $n\times n$ matrix. On $Gr(n-p,n)$ we have the universal rank $p$ quotient $Q=k^{n}\otimes \sO_{Gr}/K$, and also the rank $p$ perpendicular space $K^\perp$ to the universal subbundle $K$. Note $K^\perp \cong Q^\vee$.  It follows that
\eq{2.16}{\P(SHom(k^{n},k^{n}))^p \cong  \P(SHom(Q,Q^\vee))^p \subset \P(SHom(Q,Q^\vee)).
}
This is a fibre bundle over $Gr(n-p,n)$ with fibre $\P(Hom(k^p,k^p))^p$, the projectivized space of symmetric $p\times p$ invertible matrices. 

We can now compute $[X_{\Gamma_n}^\vee]$ as follows. Write $c(n,p) = [\P(SHom(k^n,k^n))^p]$. We have the following relations:
\ga{}{c(n,1) = [\P^{n-1}]; \quad \sum_{p=1}^n c(n,p) = [\P^{\binom{n+1}{2}-1}]; \\ c(n,p)=[Gr(n-p,n)]\cdot c(p,p) \label{2.20} \\
[X_{\Gamma_n}^\vee] = \sum_{p=1}^{n-2} c(n-1,p)
}
Here \eqref{2.20} follows from \eqref{2.16}. It is easy to see that these formulas lead to an expression for $[X_{\Gamma_n}^\vee]$ as a polynomial in the $[\P^N]$ and $[Gr(n-p-1,n-1)]$ (though the precise form of the polynomial seems complicated).  To finish the proof of the proposition, we have to show that $[Gr(a,b)] \in \Z[\A^1_k]$. Fix a splitting $k^b = k^{b-a}\oplus k^a$. Stratify $Gr(a,b) = \coprod_{p=0}^a Gr(a,b)^p$ where 
\ml{}{Gr(a,b)^p = \\
\{V \subset k^{b-a}\oplus k^a\ |\ \dim(V)=a,\ \text{Image}(V \to k^a)\text{ has rank $p$} \} = \\
\{(X,Y,f)\ |\ X \subset k^{b-a},\ Y \subset k^a,\ f:Y \to X\}
}
where $\dim X = a-p,\ \dim(Y) = p$. This is a fibration over $Gr(b-a-p,b-a)\times Gr(p,a)$ with fibre $\A^{p(b-a-p)}$. By induction, we may assume $[Gr(b-a-p,b-a)\times Gr(p,a)] \in \Z[\A^1_k]$. Since the class in the Grothendieck group of a Zariski locally trivial fibration is the class of the base times the class of the fibre, we conclude $[Gr(a,b)^p] \in \Z[\A^1_k]$, completing the proof. 
\end{proof}

In fact, we will need somewhat more. 
\begin{lem}\label{lem4.3} Let $\Gamma$ be a graph.\newline\noindent (i) Let $e_0\in \Gamma$ be an edge. Define $\Gamma' = \Gamma\cup \ve$, the graph obtained from $\Gamma$ by adding an edge $\ve$ with $\partial \ve = \partial e_0$. Then $X_{\Gamma'}^\vee$ is a cone over $X_\Gamma^\vee$. \newline\noindent (ii) Define $\Gamma' = \Gamma\cup \ve$ where $\ve$ is a tadpole, i.e. $\partial \ve = 0$. Then  $X_{\Gamma'}^\vee$ is a cone over $X_\Gamma^\vee$.
\end{lem}
\begin{proof} We prove (i). The proof of (ii) is similar and is left for the reader. 

Let $E,V$ be the edges and vertices of $\Gamma$. We have a diagram
\eq{e9}{\begin{CD} \Q^E @>\partial >> \Q^V \\
@VVV @| \\
\Q^E\oplus \Q\cdot \ve @>\partial >> \Q^V
\end{CD}
}
Dualizing and playing our usual game of interpreting edges as functionals on $\text{Image}(\partial)^\vee \cong \Q^V/\Q$, we see that $\ve^\vee = e_0^\vee$. Fix a basis for $\Q^V/\Q$ so the $(e^\vee)^2$ correspond to symmetric matrices $M_e$. We have
\eq{e10}{X_\Gamma^\vee: \det(\sum_E A_eM_e)=0;\quad X_{\Gamma'}^\vee: \det(A_\ve M_{e_0}+\sum_E A_eM_e) = 0. 
}
The second polynomial is obtained from the first by the substitution $A_{e_0} \mapsto A_{e_0}+A_\ve$. Geometrically, this is a cone as claimed. 
\end{proof}

Let $\Gamma_N$ be the complete graph on $N\ge 3$ vertices. Let $\Gamma\supset \Gamma_N$ be obtained by adding $r$ new edges (but no new vertices) to $\Gamma_N$. 
\begin{prop}\label{prop4.4} $[X_\Gamma^\vee] \in \Z[\A^1] \subset K_0(Var_k)$. 
\end{prop}
\begin{proof} Note that every pair of distinct vertices in $\Gamma_N$ are connected by an edge, so the $r$ new edges $e$ either duplicate existing edges or are tadpoles ($\partial e=0$). It follows from lemma \ref{lem4.3} that $X_\Gamma^\vee$ is an iterated cone over $\X_{\Gamma_N}^\vee$. In the Grothendieck ring, the class of a cone is the sum of the vertex point with a product of the base times an affine space, so we conclude from proposition \ref{prop4.1}. 
\end{proof}
\section{The Main Theorem}
Fix $n\ge 3$. Let $\Gamma_n$ be the complete graph on $n$ vertices. It has $\binom{n}{2}$ edges. Recall (lemma \ref{lem2.2}) a set $\{e_1,\dotsc,e_p\} \subset \text{edge}(\Gamma_n)$ is nondegenerate if cutting these edges (but leaving all vertices) does not disconnect $\Gamma_n$. (For the case $n=3$ see \eqref{2.8} and \eqref{2.9}.) Define
\eq{5.1}{S_n := \sum_{\substack{\{e_1,\dotsc,e_p\} \\
\text{nondegenerate}}} [X_{\Gamma_n-\{e_1,\dotsc,e_p\}}] \in K_0(Var_k).
}
Let $\Gamma$ be a connected graph with $n$ vertices and no multiple edges or tadpoles. Let $G \subset Sym(\text{vert}(\Gamma))$ be the subgroup of the symmetric group on the vertices which acts on the set of edges. Then $[X_\Gamma]$ appears in $S_n$ with multiplicity $n!/|G|$. 
\begin{thm} $S_n \in \Z[\A^1_k] \subset K_0(Var_k)$. 
\end{thm}
\begin{proof}It follows from \eqref{3.6} and proposition \ref{prop4.1} that 
\eq{}{\sum_{\substack{\{e_1,\dotsc,e_p\} \\
\text{nondegenerate}}} [X_{\Gamma_n-\{e_1,\dotsc,e_p\}}^0] \in \Z[\A^1_k].
}
Write $\vec{e} = \{e_1,\dotsc,e_p\}$ and let $\vec{f} = \{f_1,\dotsc,f_q\}$ be another subset of edges. We will say the pair $\{\vec{e},\vec{f}\}$ is nondegenerate if $\vec{e}$ is nondegenerate in the above sense, and if further $\vec{e}\cap \vec{f} = \emptyset$ and the edges of $\vec{f}$ do not support a loop. For $\{\vec{e},\vec{f}\}$ nondegenerate, write $(\Gamma_n-\vec{e})/\vec{f}$ for the graph obtained from $\Gamma_n$ by removing the edges in $\vec{e}$ and then contracting the edges in $\vec{f}$. If we fix a nondegenerate $\vec{e}$, we have
\eq{}{\sum_{\substack{\vec{f} \\
\{\vec{e},\vec{f}\} \text{ nondeg.}}} [X_{(\Gamma_n-\vec{e})/\vec{f}}^0] + t = [X_{\Gamma_n-\vec{e}}].
}
Here $t\in \Z[\A^1]$ accounts for the $\vec{f}$ which support a loop. These give rise to degenerate edges in $X_{\Gamma_n-\vec{e}}$ which are linear spaces and hence have classes in $ \Z[\A^1]$.  Summing now over both $\vec{e}$ and $\vec{f}$, we conclude
\eq{5.4}{S_n \equiv \sum_{\substack{\{\vec{e}, \vec{f}\} \\ \text{nondegen.}}}[X_{(\Gamma_n-\vec{e})/\vec{f}}^0] \mod \Z[\A^1]. 
}

Note that if $\vec{e}, \vec{f}$ are disjoint and $\vec{f}$ does not support a loop, then $\vec{e}$ is nondegenerate in $\Gamma_n$ if and only if it is nondegenerate in $\Gamma_n/\vec{f}$. This means we can rewrite \eqref{5.4}
\eq{5.5}{S_n \equiv \sum_{\vec{f}} \sum_{\substack{\vec{e} \subset \Gamma_n/\vec{f} \\ \text{nondegen.}}} [X_{(\Gamma_n/\vec{f})-\vec{e}}^0].
}

Let $\vec{f} = \{f_1,\dotsc,f_q\}$ and assume it does not support a loop. Then $\Gamma_n/\vec{f}$ has $n-q$ vertices, and every pair of distinct vertices is connected by at least one edge. This means we may embed $\Gamma_{n-q} \subset \Gamma_n/\vec{f}$ and think of $\Gamma_n/\vec{f}$ as obtained from $\Gamma_{n-q}$ by adding duplicate edges and tadpoles. We then apply proposition \ref{prop4.4} to conclude that $[X_{\Gamma_n/\vec{f}}^\vee] \in \Z[\A^1_k]$. Now arguing as in \eqref{3.6} we conclude
\eq{}{ \sum_{\substack{\vec{e} \subset \Gamma_n/\vec{f} \\ \text{nondegen.}}} [X_{(\Gamma_n/\vec{f})-\vec{e}}^0] \in \Z[\A^1_k]
}
Finally, plugging into \eqref{5.5} we get $S_n \in \Z[\A^1]$ as claimed. \end{proof}

\newpage \bibliographystyle{plain} \renewcommand\refname{References}

\end{document}